# $\vec{C}$ sequential optimization numbers

Zile Hui 51174500096@stu.ecnu.edu.cn

## Abstract

This work establishes a definition that is more basic than the previous ones, for the Stirling numbers of first kind, which is a sufficient but not necessary condition for the previous definition. Based on this definition and a combinatorial problem, we discover $\boldsymbol{C}$ sequential optimization numbers, where $\boldsymbol{C}$ is a $k$+1-tuple vector. For $\boldsymbol{C} = (0,1)$, we prove that $\boldsymbol{C}$ sequential optimization numbers are the unsigned Stirling numbers of first kind. We can deduce the properties of $\boldsymbol{C}$ sequential optimization numbers by following the properties of the Stirling numbers of first kind and we give specific examples such as the recurrence formula and an instance of $\boldsymbol{C}$ sequential optimization numbers. We also give specific new properties such as an explicit upper bound of them. We prove the probability that the unsigned Stirling numbers of first kind are concentrated in $O(\log n)$ is nearly 100%.

*Keywords:*

$\boldsymbol{C}$ Sequential optimization numbers, Stirling numbers of first kind, $k$-dimensional color boards problem, Explicit upper bound

## 1. Introduction

Stirling numbers were introduced by the Scottish mathematician James Stirling in his famous treatise [12]. Stirling numbers of the first kind, denoted by $s(n,m)$, are among the most important sequences in mathematics and have numerous fields such as combinatorics, number theory, numerical analysis, and probability theory [1],[3].

The recurrence formula can define the Stirling numbers of first kind. There is no known closed-form expression for the Stirling numbers of first kind and asymptotic formulas have been studied in multiple works [1][2],[6],[9],[11],[12],[19][20]. The unsigned Stirling numbers of first kind are denoted by $s_u(n,m)$. The first asymptotic formula for these numbers was given by Jordan [12]. For values of $m$ fixed as $n \to \infty$, Jordan is credited by Moser and Wyman for the asymptotic formula

$$s_u(n,m) \sim \frac{(n-1)!}{(m-1)!}(\log n + \gamma)^{m-1}$$

where $\gamma$ is Euler's constant. Moser and Wyman [15] have extended this first-order asymptotic formula to values of $m = O(\log n)$. Recently, Arratia and DeSalvo [2] gave the following results, which completely describe the asymptotic behavior of the Stirling numbers of the first kinds for various values of $n$ and $m$. For $n \to \infty$,

$$s_u(n,m) = \begin{cases} \dfrac{(n-1)!}{(m-1)!}(\log n + \gamma)^{m-1}\left(1 + O\left((\log n)^{-1}\right)\right), & m = O(\log n) \\ \dfrac{\Gamma(n+R)}{R^m \Gamma(R)\sqrt{2\pi H}}\left(1 + O(n^{-1})\right), & \sqrt{\log n} \le m \le n - n^{1/3} \\ \dbinom{n}{m}\left(\dfrac{m}{2}\right)^{n-m}\left(1 + O\left(n^{-1/3}\right)\right), & n - n^{1/3} m \le n \end{cases}$$

where $R$ is the unique solution to $\sum_{i=0}^{n-1} \frac{R}{R+i} = m$ and $H = m - \sum_{i=0}^{n-1} \frac{R^2}{(R+i)^2}$. Adell [1] gives an explicit upper bound for Stirling numbers of the first kind. For $m = 2,3,\dots,n-1$,

$$s_u(n,m) \le \frac{(n-1)!\,(\log(n-1))^{m-1}}{(m-1)!}\left(1+\frac{m-1}{\log(n-1)}\right)$$

In this work, we formally study the $k$-dimensional color boards problem and define the optimization set. We discovered the more basic definitions of the Stirling numbers of first kind. For the definition of the Stirling numbers of first kind, Theorem 4.1 gives an expression for any particular combination of $m$ elements, whereas the previous expression is for the sum of all combinations of $m$ elements. The previous definition can be derived from Theorem 4.1, but not the other way around. Based on this and further exploration of $k$-dimensional color boards problem, we discover the $\boldsymbol{C}$ sequential optimization numbers. Properties of $\boldsymbol{C}$ sequential optimization numbers can be obtained by referring to the existing properties of Stirling numbers of the first kind and we give specific examples. We study the asymptotic formulas for $\boldsymbol{C}$ sequential optimization numbers and give an explicit upper bound, which is proved by employing the recurrence formula and mathematical induction. We give the upper ratio of the sum of this upper bound to the sum of $\boldsymbol{C}$ sequential optimization numbers. We also employ this upper bound to obtain further properties, such as locating the main concentration of $\boldsymbol{C}$ sequential optimization numbers on either the head or the end for a large $n$. Particularly, for $\boldsymbol{C} = (0,1)$, the probability that the unsigned Stirling numbers of first kind are concentrated in $O(\log n)$ is nearly 100%. We transform them into each other to reveal the connections between some of the combinatorial problems in this work.

The work is organized as follows. In the following section, we study the $k$-dimensional color boards problem and define $\boldsymbol{C}$ sequential optimization numbers. Sections 3 contains several properties of $\boldsymbol{C}$ sequential optimization numbers. The more basic definition of the Stirling numbers of first kind is presented in Section 4. In Section 5, we give an explicit upper bound of $\boldsymbol{C}$ sequential optimization numbers and employ it to obtain further properties. We conclude in Section 6.

## 2. $K$-dimensional color boards problem

We start with a $k$-dimensional color boards problem in Problem 2.1. We define the optimization set in Definition 2.1 to solve this problem. Based on the optimization set, we define $\boldsymbol{C}$ sequential optimization numbers in Definition 2.2. In Theorem 2.1, we give an expression for the $\boldsymbol{C}$ sequential optimization numbers. In Lemma 2.1, we show the relationship between $\boldsymbol{C}$ sequential optimization numbers and $k$-dimensional color boards problem and give corresponding answers.

**Problem 2.1.** *As shown in Figure 1, there are $n$ boards and their order is fixed. Each board is painted one color and divided into $k$ smaller boards. For $\beta = 1,2,\dots,k$, the $\beta$th smaller board of each board form a group and their height are $1,2,\dots,n$ respectively. The problem is to determine the number of ways that $m$ colors can be seen, by following the direction of the arrow.*

We call it $k$-dimensional color boards problem. For $n = 4$ and $k = 5$, an example is shown in Figure 1a, and the number of colors that can be seen is shown in Figure 1b.

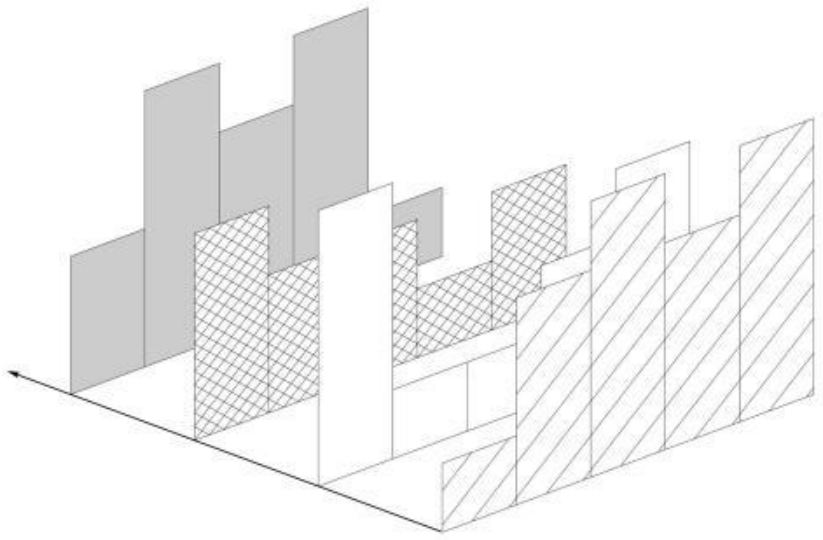

(*a*)

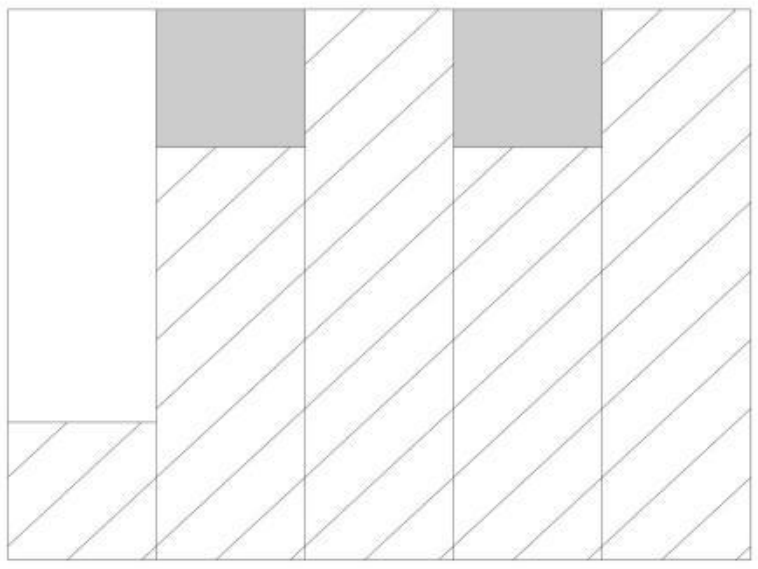

(*b*)

Figure 1: An example of *k*-dimensional color boards problem.

**Definition 2.1.** *Let* $\boldsymbol{x}(x_1, x_2, \ldots, x_k)$, $\boldsymbol{y}(y_1, y_2, \ldots, y_k)$ *be k-dimensional vectors and* $\boldsymbol{R}(R_1, R_2, \ldots, R_k)$ *be k-dimensional relation vector. If for* $\beta = 1,2,\ldots,k$, $(x_\beta, y_\beta) \in R_\beta$, *then* $\boldsymbol{x}$ *is said to be related to* $\boldsymbol{y}$ *by* $\boldsymbol{R}$, *denoted by* $(\boldsymbol{x}, \boldsymbol{y}) \in \boldsymbol{R}$. *Let U be a set of vectors,* $\boldsymbol{R}$ *be a relation vector and* $A \subseteq U$, *if for* $\forall \boldsymbol{u} \in U$, $\exists \boldsymbol{a} \in A$ *imply that* $(\boldsymbol{a}, \boldsymbol{u}) \in \boldsymbol{R}$ *or* $\boldsymbol{a} = \boldsymbol{u}$, *then A is said to be a majorization set of U by* $\boldsymbol{R}$. *If S is a majorization set of U by* $\boldsymbol{R}$ *and* $|S|$ *is a minimum of all the cardinalities of majorization sets, then S is said to be an optimization set of U by* $\boldsymbol{R}$, *denoted by* $S \overset{\boldsymbol{R}}{\subseteq} U$, *and* $|S|$ *is said to be the weight of U by* $\boldsymbol{R}$, *denoted by* $W_{\boldsymbol{R}}(U)$.

This strategy has been applied in several territories such as game theory, economics, and algorithm design [7],[8],[10],[14],[16],[17],[21]. This strategy plays a key role in multi-objective optimization [7],[10].

**Definition 2.2.** *For* $\alpha = 1,2,\ldots,n$ *and* $\beta = 1,2,\ldots,k$, $a_{1\beta}, a_{2\beta}, \ldots a_{n\beta}$ *are* $1,2,\ldots,n$ *respectively,* $\boldsymbol{a}_\alpha = (\alpha, a_{\alpha 1}, a_{\alpha 2}, \ldots a_{\alpha k})$, $\boldsymbol{b}_{\alpha\beta} = (\alpha, a_{\alpha\beta})$, $c_0$ *and* $c_\beta \in \{0,1\}$. *Let* $U = \{\boldsymbol{a}_1, \boldsymbol{a}_2, \ldots, \boldsymbol{a}_n\}$, $U_\beta = \{\boldsymbol{b}_{1\beta}, \boldsymbol{b}_{2\beta}, \ldots, \boldsymbol{b}_{n\beta}\}$, $\boldsymbol{C} = (c_0, c_1, \ldots, c_k)$ *and* $\boldsymbol{R} = (R_0, R_1, \ldots, R_k) = (<, <, \ldots, <)$ *be k+1-dimensional relation vector.* $S_\beta$ *is an optimization set of* $U_\beta$ *by* $\boldsymbol{R}_\beta(R_0, R_\beta)$. *If* $\boldsymbol{b}_{\alpha\beta} \in S_\beta$, *then* $\mu_{\alpha\beta} = 1$, *otherwise* $\mu_{\alpha\beta} = 0$. $l_\alpha = \sum_{\beta=1}^{k} \mu_{\alpha\beta}$. $S \subseteq U$. *If* $c_{l_\alpha} = 1$, *then* $\boldsymbol{a}_\alpha \in S$, *otherwise,* $\boldsymbol{a}_\alpha \notin S$. *S is said to be an* $\boldsymbol{C}$ *sequential optimization set of U by* $\boldsymbol{R}$, *denoted by* $S \overset{\boldsymbol{R},\boldsymbol{C}}{\subseteq} U$. $|S|$ *is said to be* $\boldsymbol{C}$ *sequential optimization weight of U by* $\boldsymbol{R}$, *denoted by* $W_{\boldsymbol{R},\boldsymbol{C}}(U)$. *The number of ways that* $W_{\boldsymbol{R},\boldsymbol{C}}(U) = m$ *are said to be* $\boldsymbol{C}$ *sequential optimization numbers, denoted by* $O_{\boldsymbol{C}}(n, m)$.

We provide certain formulas in preparation for later proofs. For $\beta = 0,1,\ldots,k$ and $j = 2,3,\ldots,n$, $d_{j,\beta} = \binom{k}{\beta} \frac{1}{(j-1)^\beta}$, $\boldsymbol{D}_j = (d_{j,0}, d_{j,1}, \ldots, d_{j,k})$, $\boldsymbol{C}' = (1,1,\ldots,1) - \boldsymbol{C}$ and $F_j(\boldsymbol{C}) = \boldsymbol{D}_j \boldsymbol{C}^T$. We define that $\prod_{j=m+1}^{m} F_j(\boldsymbol{C}) = 1$ and $\sum_{j=m+1}^{m} f(j) = 0$ to make the expressions concise.

(a). *For* $2 \le j_1 \le j_2 \le n$,

$$F_{j_2}(\boldsymbol{C}) \le F_{j_1}(\boldsymbol{C}) \tag{1}$$

**Proof.** For $2 \le j_1 \le j_2 \le n$,

$$\frac{d_{j_2,\beta}}{d_{j_1,\beta}} = \frac{(j_1 - 1)^\beta}{(j_2 - 1)^\beta} \le 1$$

$$d_{j_2,\beta} \le d_{j_1,\beta}$$

$$F_{j_2}(\boldsymbol{C}) \le F_{j_1}(\boldsymbol{C})$$

(b).

$$\prod_{j=2}^{n} F_j(\boldsymbol{C}) \le n^k \tag{2}$$

**Proof.**

$$\begin{aligned}\prod_{j=2}^{n} F_j(\boldsymbol{C}) &\le \prod_{j=2}^{n} F_j((1,1,\dots,1)^T)\\ &= \prod_{j=2}^{n}\sum_{\beta=0}^{k} d_{j,\beta}\\ &= \prod_{j=2}^{n}\sum_{\beta=0}^{k}\binom{k}{\beta}\frac{1}{(j-1)^\beta}\\ &= \prod_{j=2}^{n}\left(1+\frac{1}{j-1}\right)^k\\ &= n^k\end{aligned}$$

(c). *If* $c_0 = 1$,

$$F_j(\boldsymbol{C}) \ge 1 \tag{3}$$

**Proof.**

$$F_j(\boldsymbol{C}) = \boldsymbol{D}\boldsymbol{C}^T \ge c_0 d_{j,0} = 1$$

**Theorem 2.1.** *For* $k, n \in N^+$ *and* $m \in N$,

$$O_{\boldsymbol{C}}(n, m + c_k - 1) = \begin{cases} (n-1)!^k \sum\left(\prod_{i=1}^{m-1} F_{j_i}(\boldsymbol{C}) \prod_{i=1}^{n-m} F_{j_i'}(\boldsymbol{C}')\right), & 1 \le m \le n \\ 0, & otherwise \end{cases}$$

*where* $j_1, j_2, \cdots, j_{m-1}$ *are all combinations consisting of* $m-1$ *elements in set* $\{2,3,\dots,n\}$ *and* $\{j_1', j_2', \cdots, j_{n-m}'\} = \{2,3,\dots,n\} - \{j_1, j_2, \cdots, j_{m-1}\}$.

**Proof of Theorem 2.1.** Let $S$ be a set consisting of specific elements and $g(S)$ denotes the numbers of ways that $S \overset{\boldsymbol{R},\boldsymbol{C}}{\subseteq} U$. The expression for $O_{\boldsymbol{C}}(n, m)$ is initially established for the case of $\boldsymbol{C} = (0,1)$, followed by a discussion on its general applicability.

In Definition 2.2, first, for $\boldsymbol{C} = (0,1)$, we can get that $S$ is the optimization set of $U$ by $\boldsymbol{R}$. For $m = 1$, as shown in Figure 2a, place $\boldsymbol{a}_2 \sim \boldsymbol{a}_n$ and the number of ways is $(n-1)!$. So, $O_{(0,1)}(n, 1) = (n-1)!$.

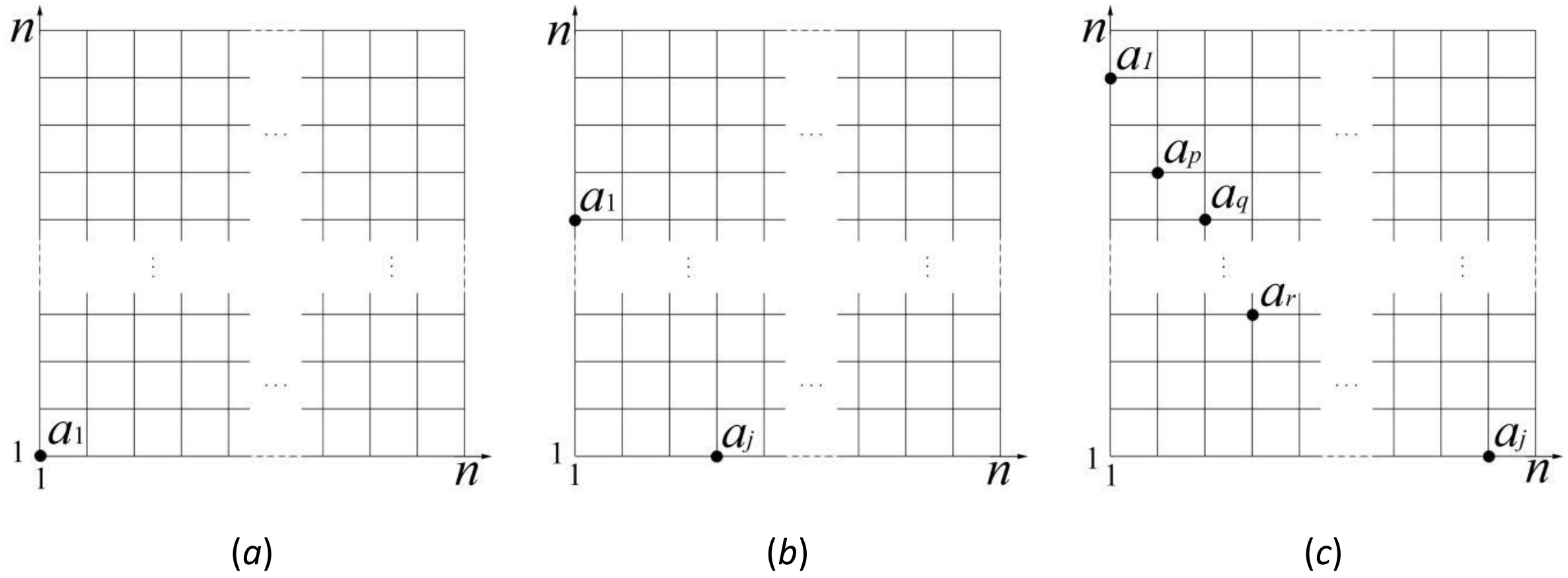

Figure 2: Three kinds of optimization set.

For $m = 2,3, \dots, n$, $1 < p < q < r < j \le n$, $S^1 = \{\boldsymbol{a}_1, \boldsymbol{a}_p, \boldsymbol{a}_q, \dots, \boldsymbol{a}_r, \boldsymbol{a}_j\}$ and $S^1$ is the optimization set of $U$ by $\boldsymbol{R}$, we apply mathematical induction to prove that $g(S^1) = \frac{(n-1)!}{(p-1)(q-1)\cdot\ldots\cdot(r-1)(j-1)}$.

Basis step: we prove $g(S^2) = \frac{(n-1)!}{j-1}$, where $S^2 = \{\boldsymbol{a}_1, \boldsymbol{a}_j\}$ and $S^2$ is optimization set of $U$ by $\boldsymbol{R}$. As shown in Figure 2b, $a_{11} = i$. We can get that $i = 2,3, \dots, n$, $j = 2,3, \dots, n$ and $i + j \le n + 2$. Place row $2 \sim (i-1)$ and the number of ways is $\frac{(n-j)!}{(n-i-j+2)!}$. Then, place row $(i+1) \sim n$ and the number of ways is $(n-i)!$. Hence,

$$
\begin{aligned}
g(S^2) &= \sum_{i=2}^{n+2-j} \frac{(n-i)!\,(n-j)!}{(n-i-j+2)!} \\
&= \frac{(n-j)!}{j-1} \sum_{i=2}^{n+2-j} \frac{(n-i)!\big((n-i+1)-(n-i-j+2)\big)}{(n-i-j+2)!} \\
&= \frac{(n-j)!}{j-1} \sum_{i=2}^{n+2-j} \left( \frac{(n-i+1)!}{(n-i-j+2)!} - \frac{(n-i)!\,(n-i-j+2)}{(n-i-j+2)!} \right)
\end{aligned}
$$

In the brackets, the right side of the minus is equal to 0 for $i = n + 2 - j$. The right side of the minus for $i = t$ is equal to the left side of the minus for $i = t + 1$, where $t = 2,3, \dots, n + 1 - j$. Hence,

$$
\begin{aligned}
g(S^2) &= \frac{(n-j)!\,(n-1)!}{(j-1)(n-j)!} \\
&= \frac{(n-1)!}{j-1}
\end{aligned}
$$

Inductive step: assume $g(S^3) = \frac{(n-1)!}{(p-1)(q-1)\cdot\ldots\cdot(r-1)}$, where $S^3 = \{\boldsymbol{a}_1, \boldsymbol{a}_p, \boldsymbol{a}_q, \dots, \boldsymbol{a}_r\}$ and $S^3$ is optimization set of $U$ by $\boldsymbol{R}$. Then, we prove $g(S^1) = \frac{(n-1)!}{(p-1)(q-1)\cdot\ldots\cdot(r-1)(j-1)}$. In Figure 2c, let $a_{r1} = i$, we can get $i = 2,3, \dots, n-1$, $j = 3,4, \dots, n$ and $i + j \le n + 2$. First, place row $2 \sim (i-1)$ and the number of ways is $\frac{(n-j)!}{(n-i-j+2)!}$. We remove all the rows and columns where the points in row $1 \sim (i-1)$ lies. In the rest of the figures, the remaining points form the set $U'$ and $S^3$ is optimization set of $U'$ by $\boldsymbol{R}$. The number of ways is $\frac{(n-i)!}{(p-1)(q-1)\cdot\ldots\cdot(r-1)}$. Hence,

$$g(S^1) = \sum_{i=2}^{n+2-j} \frac{(n-i)!\,(n-j)!}{(p-1)(q-1)\cdot\ldots\cdot(r-1)(n-i-j+2)!}$$

$$= \frac{(n-1)!}{(p-1)(q-1)\cdot\ldots\cdot(r-1)(j-1)} \tag{4}$$

From the proof of the Inductive step, we can conclude that for $\boldsymbol{C} = (0,1)$, $S^1$ can be created by adding $\boldsymbol{a}_j$ to $S^3$ and $g(S^1) = g(S^3)/(j-1)$, where $j = 2,3,\ldots,n$.

Second, for $k > 1$, let $S^1 \overset{\boldsymbol{R},\boldsymbol{C}}{\subseteq} U$ and $S^3 \overset{\boldsymbol{R},\boldsymbol{C}}{\subseteq} U$. We discuss the relationship between $g(S^1)$ and $g(S^3)$. We find that the relationship between $S_\beta$ and $g(S_\beta)$ is the same as the relationship between $S$ and $g(S)$ for $\boldsymbol{C} = (0,1)$. Let $S^{1,3} = \{\boldsymbol{a}_1, \boldsymbol{a}_p, \boldsymbol{a}_q, \ldots, \boldsymbol{a}_r\}$ if $c_0 = 0$, otherwise $S^{1,3} = \{\boldsymbol{a}_1, \boldsymbol{a}_p, \boldsymbol{a}_q, \ldots, \boldsymbol{a}_r, \boldsymbol{a}_j\}$. $S^{1,3} \overset{\boldsymbol{R},\boldsymbol{C}}{\subseteq} U$ and there are no $\boldsymbol{b}_{j\beta}$ in $S_\beta^{1,3}$. We can create $l$ $S_\beta^1$ by adding $\boldsymbol{b}_{j\beta}$ to $S_\beta^{1,3}$, which can create $S^1$ by $S^{1,3}$ for $c_l = 1$. We can create $l$ $S_\beta^3$ by adding $\boldsymbol{b}_{j\beta}$ to $S_\beta^{1,3}$, which can create $S^3$ by $S^{1,3}$ for $c_l = 0$. Hence,

$$g(S^1) = \boldsymbol{D}\boldsymbol{C}^T g(S^{1,3})$$

$$g(S^3) = \boldsymbol{D}\boldsymbol{C}'^T g(S^{1,3})$$

and

$$g(S^1) = \frac{\boldsymbol{D}\boldsymbol{C}^T}{\boldsymbol{D}\boldsymbol{C}'^T} g(S^3) = \frac{F_j(\boldsymbol{C})}{F_j(\boldsymbol{C}')} g(S^3) \tag{5}$$

We mention $c_k$ in two cases since $\boldsymbol{b}_{1\beta}$ must be in $S_\beta$.

(a). For $c_k = 1$,

$$O_{\boldsymbol{C}}(n,m) = \begin{cases} (n-1)!^k \sum \prod_{i=1}^{n-1} F_{j_i'}(\boldsymbol{C}'), & m = 1 \\ (n-1)!^k \sum \left( \prod_{i=1}^{m-1} F_{j_i}(\boldsymbol{C}) \prod_{i=1}^{n-m} F_{j_i'}(\boldsymbol{C}') \right), & 2 \le m \le n-1 \\ (n-1)!^k \sum \prod_{i=1}^{n-1} F_{j_i}(\boldsymbol{C}), & m = n \\ 0, & m = 0 \text{ or } m > n \end{cases}$$

where $j_1, j_2, \cdots, j_{m-1}$ are all combinations consisting of $m-1$ elements in set $\{2,3,\ldots,n\}$ and $\{j_1', j_2', \cdots, j_{n-m}'\} = \{2,3,\ldots,n\} - \{j_1, j_2, \cdots, j_{m-1}\}$.

(b). For $c_k = 0$,

$$O_{\boldsymbol{C}}(n,m)=\begin{cases}(n-1)!^k\sum\prod_{i=1}^{n-1}F_{j_i'}(\boldsymbol{C}'), & m=0\\(n-1)!^k\sum\left(\prod_{i=1}^{m}F_{j_i}(\boldsymbol{C})\prod_{i=1}^{n-m-1}F_{j_i'}(\boldsymbol{C}')\right), & 1\le m\le n-2\\(n-1)!^k\sum\prod_{i=1}^{n-1}F_{j_i}(\boldsymbol{C}), & m=n-1\\0, & m=-1\text{ or }m>n-1\end{cases}$$

where $j_1, j_2, \cdots, j_m$ are all combinations consisting of $m$ elements in set $\{2,3,\ldots,n\}$ and $\{j_1', j_2', \cdots, j_{n-m-1}'\}=\{2,3,\ldots,n\}-\{j_1, j_2, \cdots, j_m\}$.

The above formula also holds for $\boldsymbol{C}=(0,0)$, $(1,0)$ and $(1,1)$. To sum up, for $k,n\in N^+$ and $m\in N$,

$$O_{\boldsymbol{C}}(n,m+c_k-1)=\begin{cases}(n-1)!^k\sum\left(\prod_{i=1}^{m-1}F_{j_i}(\boldsymbol{C})\prod_{i=1}^{n-m}F_{j_i'}(\boldsymbol{C}')\right), & 1\le m\le n\\0, & \text{otherwise}\end{cases}$$

where $j_1, j_2, \cdots, j_{m-1}$ are all combinations consisting of $m-1$ elements in set $\{2,3,\ldots,n\}$ and $\{j_1', j_2', \cdots, j_{n-m}'\}=\{2,3,\ldots,n\}-\{j_1, j_2, \cdots, j_{m-1}\}$.

**Lemma 2.1.** *In Problem 2.1, following the direction of the arrow, the number of ways that m colors can be seen is* $O_{(0,1,1,\ldots,1)}(n,m)$. *For* $\boldsymbol{C}=(c_0,c_1,c_2,\ldots,c_k)$, $c_{v_1}, c_{v_2},\ldots,c_{v_w}$ *are elements which are equal to 1 in* $\boldsymbol{C}$. *Following the direction of the arrow, the number of ways that m colors can be seen* $v_1$ *or* $v_2$ *or…or* $v_w$ *times is* $O_{\boldsymbol{C}}(n,m)$.

**Proof of lemma 2.1.** For $\alpha=1,2,\ldots,n$ and $\beta=1,2,\ldots,k$, label as the $\alpha$th board follow direction of the arrow. Label the height of the smaller boards in same color from left to right as $h_{\alpha1}, h_{\alpha2},\ldots h_{\alpha k}$. Let $\boldsymbol{R}=(<,>,>,\ldots,>)$ be *k+1*-dimensional relation vector, $\boldsymbol{h}_\alpha=(\alpha,h_{\alpha1},h_{\alpha2},\ldots h_{\alpha k})$ and $a_{\alpha\beta}=n+1-h_{\alpha\beta}$. The following proof refers to the proof of Theorem 2.1. For $\boldsymbol{C}=(0,1)$, as shown in Figure 3, we can get same optimization set in Figure 2c follow direction of the line of sight. For $k>1$, $j$th color is seen $l$ times implies that we add $l$ $\boldsymbol{b}_{j\beta}$ to $S_\beta^{1,3}$. First color is seen $k$ times and $\boldsymbol{b}_{1\beta}$ must be in $S_\beta$. Hence, the number of ways that $m$ colors can be seen $v_1$ or $v_2$ or … or $v_w$ times is $O_{\boldsymbol{C}}(n,m)$ and this problem can be considered as an instance of $\boldsymbol{C}$ sequential optimization numbers. One color being seen is equal to that of it being seen 1 or 2 or … or $k$ times. Hence, we can conclude that the number of ways that $m$ colors can be seen is $O_{(0,1,1,\ldots,1)}(n,m)$.

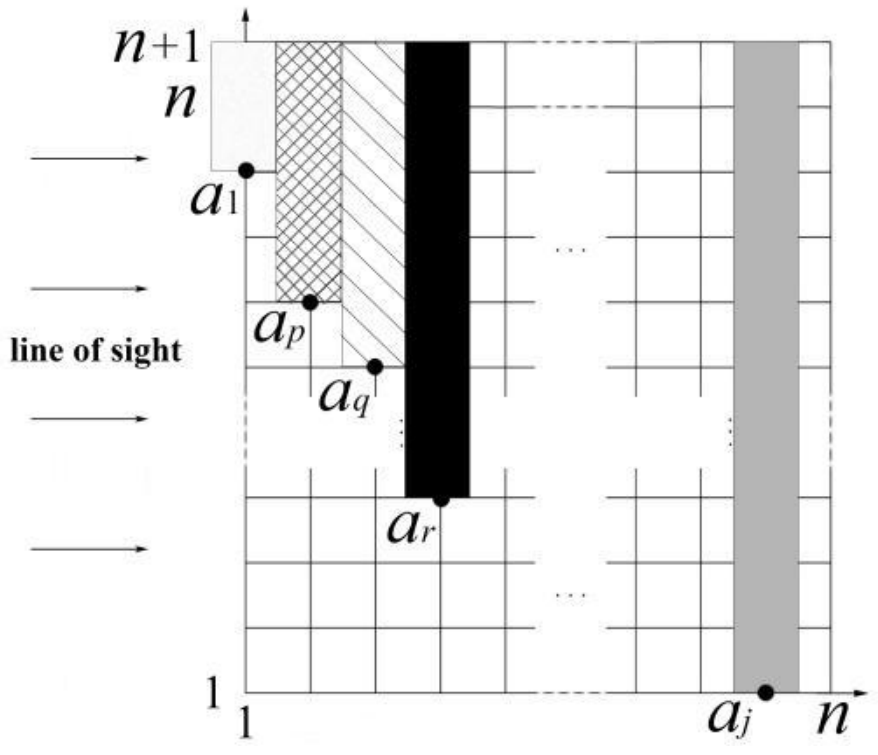

Figure 3: Relationship between smaller color boards group and optimization set.

The combinatorics community has accumulated the knowledge of multiple number sequences summarized in the On-Line Encyclopedia of integer Sequences (OEIS). We find that 2-dimensional color boards problem has been presented as sequence A309053 in OEIS. In sequence A309053, Duff calculates this sequence by enumeration and the computational complexity is at least of the order of exponential. However, we can apply the recurrence formula in Theorem 3.1 to calculate the same, where the computational complexity reduces to that of a polynomial. Apart from the unsigned Stirling numbers of first kind and sequence A309053, we find no other sequences in $\boldsymbol{C}$ Sequential optimization numbers in OEIS. This suggests that $\boldsymbol{C}$ Sequential optimization numbers have not been studied previously.

## 3. Several properties of $C$ sequential optimization numbers

Lemma 3.1 shows that there is a duality inside $\boldsymbol{C}$ sequential optimization numbers. Lemma 3.2 is a special case of $O_{\boldsymbol{C}}(n,m)$. Lemma 3.3 shows that the unsigned Stirling numbers of first kind are (0,1) sequential optimization numbers. Thereafter, we give certain properties of $\boldsymbol{C}$ sequential optimization numbers by referring to the common properties of the Stirling numbers of first kind. In Theorem 3.1, we give the recurrence formula for $\boldsymbol{C}$ sequential optimization numbers. For $\boldsymbol{C} = (0,1)$, Lemma 3.4 usually appears in the introduction of the Stirling numbers of first kind.

**Lemma 3.1.** *For* $c_k - 1 \le m \le n + c_k$,

$$O_{\boldsymbol{C}}(n,m) = O_{\boldsymbol{C}'}(n, n-m)$$

*For* $0 \le m \le n+1$,

$$O_{\boldsymbol{C}}(n, m + c_k - 1) = O_{\boldsymbol{C}'}(n, n - m + c_k').$$

**Proof of Lemma 3.1**. For $1 \le t \le n$,

$$\begin{aligned} O_{\boldsymbol{C}}(n, t + c_k - 1) &= (n-1)!^k \sum \left( \prod_{i=1}^{t-1} F_{j_i}(\boldsymbol{C}) \prod_{i=1}^{n-t} F_{j_i'}(\boldsymbol{C}') \right) \\ &= (n-1)!^k \sum \left( \prod_{i=1}^{n-t+1-1} F_{j_i'}(\boldsymbol{C}') \prod_{i=1}^{n-(n-t+1)} F_{j_i}(\boldsymbol{C}) \right) \\ &= O_{\boldsymbol{C}'}(n, n - t + 1 + c_k' - 1) \\ &= O_{\boldsymbol{C}'}(n, n - t + c_k') \end{aligned}$$

For $t = 0$ and $t = n+1$,

$$O_{\boldsymbol{C}}(n, t + c_k - 1) = O_{\boldsymbol{C}'}(n, n - t + c_k') = 0$$

Let $m = t + c_k - 1$,

$$O_{\boldsymbol{C}}(n,m) = O_{\boldsymbol{C}'}(n, n-m)$$

where $c_k - 1 \le m \le n + c_k$.

**Lemma 3.2.**

$$O_{(0,0,\dots,0)}(n,0) = O_{(1,1,\dots,1)}(n,n) = \sum_{m=0}^{n} O_{\boldsymbol{C}}(n,m) = n!^k$$

**Proof of Lemma 3.2.** In Theorem 2.1, $O_{(0,0,\dots,0)}(n,0) = O_{(1,1,\dots,1)}(n,n) = n!^k$. For $m < 0$ and $m > n$, $O_{\boldsymbol{C}}(n,m) = 0$. So, $\sum_{m=0}^{n} O_{\boldsymbol{C}}(n,m)$ is the number of permutations.

**Lemma 3.3.** *For* $\boldsymbol{C} = (0,1)$, $O_{\boldsymbol{C}}(n,m) = s_u(n,m)$, *where* $s_u(n,m)$ *are the unsigned Stirling numbers of first kind.*

**Proof of Lemma 3.3** In Theorem 2.1, let $\boldsymbol{C} = (0,1)$ and we can get

$$O_{(0,1)}(n,m) = \begin{cases} (n-1)! \sum \prod_{i=1}^{m-1} \dfrac{1}{j_i - 1}, & 1 \le m \le n \\ 0, & \text{otherwise} \end{cases}$$

where $j_1, j_2, \cdots, j_{m-1}$ are all combinations consisting of $m-1$ elements in set $\{2,3,\dots,n\}$. This is a definition of the unsigned Stirling numbers of first kind.

**Theorem 3.1.** *For* $k, n \in N^+$ *and* $m \ge c_k - 1$, *the recurrence formula of* $\boldsymbol{C}$ *sequential optimization numbers is*

$$O_{\boldsymbol{C}}(n+1, m+1) = n^k F_{n+1}(\boldsymbol{C}) O_{\boldsymbol{C}}(n,m) + n^k F_{n+1}(\boldsymbol{C}') O_{\boldsymbol{C}}(n, m+1)$$

*and the boundary condition is*

$$O_c(n,m) = \begin{cases} 0, & m = c_k - 1 \text{ or } m > c_k - 1 + n \\ 1, & m = c_k, n = 1 \end{cases}$$

**Proof of Theorem 3.1.** For $c_k = 1, n \ge 2$ and $m = 1,2,\dots,n-1$. we divide the number of ways that $O_k(n+1, m+1)$ denotes into two parts, with and without $\boldsymbol{a}_{n+1}$.

First, ways with $\boldsymbol{a}_{n+1}$ and $m$ vectors from $\{\boldsymbol{a}_1, \boldsymbol{a}_2, \dots, \boldsymbol{a}_n\}$. Let $j_m$ be $n+1$, $j_1, j_2, \cdots, j_{m-1}$ be all combinations consisting of $m-1$ elements in set $\{2,3,\dots,n\}$ and $\{j'_1, j'_2, \cdots, j'_{n-m}\} = \{2,3,\dots,n\} - \{j_1, j_2, \cdots, j_{m-1}\}$.

$$\begin{aligned} n!^k \sum \left( \prod_{i=1}^{m} F_{j_i}(\boldsymbol{C}) \prod_{i=1}^{n-m} F_{j'_i}(\boldsymbol{C}') \right) &= n!^k\, F_{n+1}(\boldsymbol{C}) \sum \left( \prod_{i=1}^{m-1} F_{j_i}(\boldsymbol{C}) \prod_{i=1}^{n-m} F_{j'_i}(\boldsymbol{C}') \right) \\ &= n^k F_{n+1}(\boldsymbol{C}) (n-1)!^k \sum \left( \prod_{i=1}^{m-1} F_{j_i}(\boldsymbol{C}) \prod_{i=1}^{n-m} F_{j'_i}(\boldsymbol{C}') \right) \\ &= n^k F_{n+1}(\boldsymbol{C}) O_{\boldsymbol{C}}(n,m) \end{aligned}$$

Then, the ways without $\boldsymbol{a}_{n+1}$ and $m$ vectors from $\{\boldsymbol{a}_1, \boldsymbol{a}_2, \dots, \boldsymbol{a}_n\}$. let $j'_{n-m}$ be $n+1$, $j_1, j_2, \cdots, j_m$ are all combinations consisting of $m$ elements in set $\{2,3,\dots,n\}$and $\{j'_1, j'_2, \cdots, j'_{n-m-1}\} = \{2,3,\dots,n\} - \{j_1, j_2, \cdots, j_m\}$.

$$\begin{aligned} n!^k \sum \left( \prod_{i=1}^{m} F_{j_i}(\boldsymbol{C}) \prod_{i=1}^{n-m} F_{j'_i}(\boldsymbol{C}') \right) &= n!^k\, F_{n+1}(\boldsymbol{C}') \sum \left( \prod_{i=1}^{m} F_{j_i}(\boldsymbol{C}) \prod_{i=1}^{n-m-1} F_{j'_i}(\boldsymbol{C}') \right) \\ &= n^k F_{n+1}(\boldsymbol{C}') (n-1)!^k \sum \left( \prod_{i=1}^{m} F_{j_i}(\boldsymbol{C}) \prod_{i=1}^{n-m-1} F_{j'_i}(\boldsymbol{C}') \right) \\ &= n^k F_{n+1}(\boldsymbol{C}') O_{\boldsymbol{C}}(n, m+1) \end{aligned}$$

To sum up,

$$O_{\boldsymbol{C}}(n+1, m+1) = n^k F_{n+1}(\boldsymbol{C}) O_{\boldsymbol{C}}(n,m) + n^k F_{n+1}(\boldsymbol{C}') O_{\boldsymbol{C}}(n, m+1)$$

According to Theorem 2.1, the boundary condition is

$$O_{\boldsymbol{c}}(n,m)=\begin{cases}0, & m=c_k-1 \text{ or } m>c_k-1+n\\ 1, & m=c_k, n=1\end{cases}$$

The same method can be iterated for $c_k=0$. The above formula also holds for other cases. Hence, the recurrence formula works for $k,n\in N^+$ and $m\geq c_k-1$.

**Lemma 3.4.** *For $k\in N^+$ and $n\geq 2$, we define $x_{\boldsymbol{C}}^{1\uparrow}=x$, $O_{\boldsymbol{C}}^{u}(n,m)=O_{\boldsymbol{C}}(n,m+c_k-1)$ and*

$$x_{\boldsymbol{C}}^{n\uparrow}=x\big(F_2(\boldsymbol{C})x+F_2(\boldsymbol{C}')\big)\left(2^kF_3(\boldsymbol{C})x+2^kF_3(\boldsymbol{C}')\right)\cdots\left((n-1)^kF_n(\boldsymbol{C})x+(n-1)^kF_n(\boldsymbol{C}')\right)$$

*We can get*

$$x_{\boldsymbol{C}}^{n\uparrow}=\sum_{m=0}^{n}O_{\boldsymbol{C}}^{u}(n,m)x^m$$

*and the zero points are $x=0$ and $x=-F_m(\boldsymbol{C}')/F_m(\boldsymbol{C})$, where $m=2,3,,\dots,n$. We call $O_{\boldsymbol{C}}^{u}(n,m)$ unsigned $\boldsymbol{C}$ sequential optimization numbers.*

*For $k\in N^+$ and $n\geq 2$, we define $x_{\boldsymbol{C}}^{1\downarrow}=x$, $O_{\boldsymbol{C}}^{s}(n,m)=(-1)^{n+m}O_{\boldsymbol{C}}(n,m+c_k-1)$ and*

$$x_{\boldsymbol{C}}^{n\downarrow}=x\big(F_2(\boldsymbol{C})x-F_2(\boldsymbol{C}')\big)\left(2^kF_3(\boldsymbol{C})x-2^kF_3(\boldsymbol{C}')\right)\cdots\left((n-1)^kF_n(\boldsymbol{C})x-(n-1)^kF_n(\boldsymbol{C}')\right)$$

*We can get*

$$x_{\boldsymbol{C}}^{n\downarrow}=\sum_{m=0}^{n}O_{\boldsymbol{C}}^{s}(n,m)x^m$$

*and the zero points are $x=0$ and $x=F_m(\boldsymbol{C}')/F_m(\boldsymbol{C})$, where $m=2,3,,\dots,n$. We call $O_{\boldsymbol{C}}^{s}(n,m)$ signed $\boldsymbol{C}$ sequential optimization numbers.*

**Proof of Lemma 3.4.** We prove it by applying Theorem 3.1.

## 4. More basic definition

In this section, we show that the more basic definitions of the Stirling numbers of first kind and $\boldsymbol{C}$ sequential optimization numbers.

**Theorem 4.1.** *Let $\boldsymbol{a}_1(1,a_1)$, $\boldsymbol{a}_2(2,a_2)$, …, $\boldsymbol{a}_n(n,a_n)$ be n 2-dimensional vectors and $a_1,a_2,\dots a_n$ be 1,2, …, $n$, respectively. $\boldsymbol{R}=(<,<)$. $U=\{\boldsymbol{a}_1,\boldsymbol{a}_2,\dots,\boldsymbol{a}_n\}$. $S$ is an optimization set of $U$ by $\boldsymbol{R}$. $g(S)$ denotes the numbers of ways, where $S$ is a set consisting of specific elements. For $S=\{\boldsymbol{a}_1,\boldsymbol{a}_p,\boldsymbol{a}_q,\dots,\boldsymbol{a}_r\}$, we can deduce the following formula from formula 4 in the Proof of Theorem 2.1,*

$$g(S)=\frac{(n-1)!}{(p-1)(q-1)\cdot\ldots\cdot(r-1)}$$

*denoted by*

$$s_n(p,q,\dots,r)=\frac{(n-1)!}{(p-1)(q-1)\cdot\ldots\cdot(r-1)}$$

For $\boldsymbol{C}=(0,1)$, Theorems 2.1 and 3.1 are the two definitions of the unsigned Stirling numbers of first kind and these two definitions are mutually sufficient and necessary [5],[12]. The following is a definition of the unsigned Stirling numbers of first kind.

$$s_u(n,m) = \begin{cases} (n-1)! \sum \prod_{i=1}^{m-1} \frac{1}{j_i - 1}, & m = 1,2,\ldots n \\ 0, & otherwise \end{cases}$$

where $j_1, j_2, \cdots, j_{m-1}$ are all combinations consisting of $m-1$ elements in set $\{2,3,\ldots,n\}$.

Theorem 4.1 gives an expression for any particular combination of $m$ elements, whereas the previous definition is an expression for the sum of all combinations of $m$ elements. The previous definition can be derived from Theorem 4.1, but not the other way around. Theorem 4.1 demonstrates that our proof entails a more basic definition than the previous proof and our definition is a sufficient but not necessary condition for the previous definition. This definition can help in discovering new sequences and studying the properties of the Stirling numbers of first kind and $\boldsymbol{C}$ sequential optimization numbers is an example.

We employ the following example to show that the more basic definition of the Stirling numbers of first kind is hidden in a problem. Circular permutations problem generally serves as an introduction to the unsigned Stirling numbers of first kind [2],[4],[5]. Further, $s_u(n,m)$ counts the number of arrangements of $n$ objects into $m$ nonempty circular permutations. We include certain restrictions to Circular permutations problem so that this problem can be equivalent to 1-dimensional color boards problem.

Label $n$ objects as $1,2,\cdots,n$ respectively and the largest number of $m$ circular permutations as $t_1, t_2, \cdots, t_m$. Disconnect $m$ circular permutations into $m$ permutations at each maximum number of them. Label the $m$ permutations objects as $v_1, v_2, \cdots, v_m$. For $i = 1,2,\cdots,m$, $t_i$ is at the front of the $v_i$. We arrange $v_1, v_2, \cdots, v_m$ into a permutation in the ascending order from $t_1, t_2, \cdots, t_m$, denoted by $A_n$. We can obtain the positions of $t_1, t_2, \cdots, t_m$ in $A_n$, denoted by $w_1, w_2, \cdots, w_m$. The number of ways that the position of $t_1, t_2, \cdots, t_m$ in $A_n$ is $w_1, w_2, \cdots, w_m$ is $s_n(w_1, w_2, \ldots, w_m)$.

**Proof.** Referring to Figure *3*, we expand $A_n$ from left to right, corresponding to $n$ different colored boards. The number in $A_n$ is the height of the board. $t_1$ is on the far left. For $i = 2,3,\cdots,m$, $t_i$ is greater than $t_{i-1}$ and $t_{i-1}$ is greater than all the numbers in $v_{i-1}$. Hence, for $i = 1,2,\cdots,m$, $t_i$ can be seen, and the remaining numbers in $v_i$ cannot be seen since they are smaller than $t_i$ and arranged to the right of $t_i$. We can get that only $t_1, t_2, \cdots, t_m$ can be seen and they are in positions $w_1, w_2, \cdots, w_m$. Hence, the number of ways that the position of $t_1, t_2, \cdots, t_m$ in $A_n$ is $w_1, w_2, \ldots, w_m$ is $s_n(w_1, w_2, \ldots, w_m)$.

Those restrictions follow from an experiment that transforms the Circular permutations problem into 1-dimensional color boards problem. In the absence of a more basic definition for the 1-dimensional color boards problem, we would not design this experiment and get the restrictions.

According to the example above, the restrictions are sophisticated. Multiple simple problems hide this more basic definition, which may be one of the reasons why this more basic definition has not been found so far.

The definition, which is a sufficient but not necessary condition for the previous definition, also applies to $\boldsymbol{C}$ sequential optimization numbers. We denote

$$O_{\boldsymbol{C},n}(p,q,\ldots,r) = (n-1)!^k \frac{F_p(\boldsymbol{C})F_q(\boldsymbol{C})\ldots F_r(\boldsymbol{C})}{F_p(\boldsymbol{C}')F_q(\boldsymbol{C}')\ldots F_r(\boldsymbol{C}')} \prod_{i=2}^{n} F_i(\boldsymbol{C}')$$

In Definition 2.2, let $S = \{\boldsymbol{a_1}, \boldsymbol{a_p}, \boldsymbol{a_q}, \ldots, \boldsymbol{a_r}\}$ if $c_k = 1$, otherwise $S = \{\boldsymbol{a_p}, \boldsymbol{a_q}, \ldots, \boldsymbol{a_r}\}$; $g(S)$ denotes the numbers of ways, where $S$ is a set consisting of specific elements. We can deduce $g(S) = O_{\boldsymbol{C},n}(p,q,\ldots,r)$ from formula 5 in the Proof of Theorem 2.1.

## 5. Explicit upper bound

In this section, we give an explicit upper bound of $\boldsymbol{C}$ sequential optimization numbers and employ it to obtain further properties. In Theorem 5.1, we give the upper bound of $\boldsymbol{C}$ sequential optimization numbers. In Lemma 5.1, for $C = (0,1)$, we give the upper bound of the unsigned Stirling numbers of first kind. Theorem 5.2 shows that the $\boldsymbol{C}$ sequential optimization numbers are nearly concentrated in the head or the end for a large $n$. Lemma 5.2 shows the probability that the unsigned Stirling numbers of first kind are concentrated in $O(\log n)$ is nearly 100%. In Theorem 5.3, we obtain the upper ratio of the sum of the upper bound to the sum of $\boldsymbol{C}$ sequential optimization numbers.

**Theorem 5.1.** *For $k \geq 1$, $n \geq 2$ and $1 \leq m \leq n$, we define*

$$O_{\boldsymbol{C}max}(n, m + c_k - 1) = \frac{(n-1)!^k}{(m-1)!}(\boldsymbol{H}_n\boldsymbol{C}^T)^{m-1}\prod_{i=1}^{n-m}F_{i+1}(\boldsymbol{C}')$$

*where $\boldsymbol{H}_n = (h_0, h_1, \dots, h_k)$ and for $\beta = 0,1, \dots, k$, $h_\beta = \binom{k}{\beta}\sum_{j=1}^{n-1}\frac{1}{j^\beta}$. We can get $O_{\boldsymbol{C}max}(n, m + c_k - 1) \geq O_{\boldsymbol{C}}(n, m + c_k - 1)$.*

**Proof of Theorem 5.1.** We can get

$$h_\beta = \binom{k}{\beta}\sum_{j=1}^{n-1}\frac{1}{j^\beta} = \sum_{j=2}^{n}\binom{k}{\beta}\frac{1}{(j-1)^\beta} = \sum_{j=2}^{n} d_{j,\beta}$$

$$\boldsymbol{H}_n\boldsymbol{C}^T = \sum_{\beta=1}^{k} h_\beta c_\beta = \sum_{\beta=1}^{k}\sum_{j=2}^{n} d_{j,\beta}c_\beta = \sum_{j=2}^{n}\sum_{\beta=1}^{k} d_{j,\beta}c_\beta = \sum_{j=2}^{n} F_j(\boldsymbol{C}) \tag{6}$$

We apply mathematical induction to prove it. First, we prove the boundary condition. For $k \geq 1$ and $n \geq 2$,

$$O_{\boldsymbol{C}max}(n, \mathrm{m} + c_k - 1)\begin{cases} = (n-1)!^k \prod_{i=1}^{n-1} F_{i+1}(\boldsymbol{C}'), & m = 1 \\ > 0, & m = n+1 \end{cases}$$

Hence, the boundary condition satisfies $O_{\boldsymbol{C}max}(n, m + c_k - 1) \geq O_{\boldsymbol{C}}(n, m + c_k - 1)$.

Then, we prove the inequality for $2 \leq m \leq n$.

Basis step: For $n = 2$,

$$O_{\boldsymbol{C}max}(2,2 + c_k - 1) = \boldsymbol{H}_2\boldsymbol{C}^T = F_2(\boldsymbol{C}) = O_{\boldsymbol{C}}(2,2 + c_k - 1) \qquad formula\ (6)$$

So, for $n = 2$, $O_{\boldsymbol{C}max}(n, m + c_k - 1) \geq O_{\boldsymbol{C}}(n, m + c_k - 1)$.

Inductive step: For $n \geq 2$ and $2 \leq m \leq n$, we assume $O_{\boldsymbol{C}max}(n, m + c_k - 1) \geq O_{\boldsymbol{C}}(n, m + c_k - 1)$. Then, we prove $O_{\boldsymbol{C}max}(n + 1, m + c_k - 1) \geq O_{\boldsymbol{C}}(n + 1, m + c_k - 1)$.

$$O_{\boldsymbol{C}max}(n + 1, \mathrm{m} + c_k - 1)$$

$$= \frac{n!^k}{(m-1)!}(\boldsymbol{H}_{n+1}\boldsymbol{C}^T)^{m-1}\prod_{i=1}^{n+1-\mathrm{m}} F_{i+1}(\boldsymbol{C}')$$

$$= \frac{n!^k}{(m-1)!}\left(\boldsymbol{H}_n\boldsymbol{C}^T + F_{n+1}(\boldsymbol{C})\right)^{m-1}\prod_{i=1}^{n+1-\mathrm{m}} F_{i+1}(\boldsymbol{C}') \qquad formula\ (6)$$

$$\geq \frac{n!^k}{(m-1)!}(\boldsymbol{H}_n\boldsymbol{C}^T)^{m-1}\prod_{i=1}^{n+1-\mathrm{m}} F_{i+1}(\boldsymbol{C}') + \frac{n!^k\,(m-1)}{(m-1)!}F_{n+1}(\boldsymbol{C})(\boldsymbol{H}_n\boldsymbol{C}^T)^{m-2}\prod_{i=1}^{n+1-\mathrm{m}} F_{i+1}(\boldsymbol{C}')$$

$$= F_{n+2-\mathrm{m}}(\boldsymbol{C}')\frac{n^k(n-1)!^k}{(m-1)!}(\boldsymbol{H}_n\boldsymbol{C}^T)^{m-1}\prod_{i=1}^{n-\mathrm{m}} F_{i+1}(\boldsymbol{C}') + \frac{n^k(n-1)!^k}{(m-2)!}F_{n+1}(\boldsymbol{C})(\boldsymbol{H}_n\boldsymbol{C}^T)^{m-2}\prod_{i=1}^{n+1-\mathrm{m}} F_{i+1}(\boldsymbol{C}')$$

$$\geq F_{n+1}(\boldsymbol{C}')\frac{n^k(n-1)!^k}{(m-1)!}(\boldsymbol{H}_n\boldsymbol{C}^T)^{m-1}\prod_{i=1}^{n-\mathrm{m}} F_{i+1}(\boldsymbol{C}') + \frac{n^k(n-1)!^k}{(m-2)!}F_{n+1}(\boldsymbol{C})(\boldsymbol{H}_n\boldsymbol{C}^T)^{m-2}\prod_{i=1}^{n+1-\mathrm{m}} F_{i+1}(\boldsymbol{C}') \qquad formula\ (1)$$

$$= n^kF_{n+1}(\boldsymbol{C}')O_{\boldsymbol{C}max}(n,m+c_k-1) + n^kF_{n+1}(\boldsymbol{C})O_{\boldsymbol{C}max}(n,m-1+c_k-1)$$

$$\geq n^kF_{n+1}(\boldsymbol{C}')O_{\boldsymbol{C}}(n,m+c_k-1) + n^kF_{n+1}(\boldsymbol{C})O_{\boldsymbol{C}}(n,m-1+c_k-1)$$

$$= O_{\boldsymbol{C}}(n+1,\mathrm{m}+c_k-1) \qquad Theorem\ \mathbf{3.1}$$

We can do same thing for $m = n+1$. To sum up, for $k \geq 1, n \geq 2$ and $1 \leq m \leq n$, $O_{\boldsymbol{C}max}(n,m+c_k-1) \geq O_{\boldsymbol{C}}(n,m+c_k-1)$.

**Lemma 5.1.** *In Theorem 5.1, let $C = (0,1)$, for $n \geq 2$ and $1 \leq m \leq n$, we can get*

$$s_u(n,m) \leq \frac{(n-1)!}{(m-1)!}\big(H(n-1)\big)^{m-1}$$

*where $H(n) = \sum_{i=1}^{n}\frac{1}{i}$.*

This formula is similar to the ones given in Section 1 and better in certain ways. Nevertheless, the proof of this formula is simple and we only use the recurrence formula of the unsigned Stirling numbers of first kind, mathematical induction, and binomial expansion. The proof is crude and leaves much room for improvement. It is possible to obtain better results of asymptotic formulas for the Stirling numbers of first kind based on this method. The expression of the upper bound leads to further results such as the following properties.

**Theorem 5.2.** *Let $n \geq 2$, $k \geq 1$, $M_1 \in N^+$ and $P\big(O_{\boldsymbol{C}}(n,m+c_k-1)\big) = \frac{O_{\boldsymbol{C}}(n,m+c_k-1)}{n!^k}$ be the probability of $O_{\boldsymbol{C}}(n,m+c_k-1)$.*

*For $c_0 = 0$ and $M = \left\lceil ekc_1(log(n-1)+\gamma) + \frac{e\pi^2}{6}\sum_{\beta=2}^{k} c_\beta \binom{k}{\beta}\right\rceil + M_1$, we can get*

$$P\big(O_{\boldsymbol{C}}(n,m > M)\big) \leq exp(-M_1)$$

*For $c_0 = 1$ and $M' = \left\lceil ekc_1'(log(n-1)+\gamma) + \frac{e\pi^2}{6}\sum_{\beta=2}^{k} c_\beta' \binom{k}{\beta}\right\rceil + M_1$, we can get*

$$P\big(O_{\boldsymbol{C}}(n,m < n+1-M')\big) \leq exp(-M_1)$$

**Proof of Theorem 5.2.** Let $P\big(O_{\boldsymbol{C}max}(n,m+c_k-1)\big) = \frac{O_{\boldsymbol{C}max}(n,m+c_k-1)}{n!^k}$ be the probability of the upper bound of $\boldsymbol{C}$ sequential optimization numbers. For $c_0 = 0$ and $1 \leq m \leq n-1$,

$$
\begin{aligned}
P\big(O_{\boldsymbol{C}max}(n, m + c_k)\big) &= \frac{O_{\boldsymbol{C}max}(n, m + c_k)}{n!^k} \\
&= \frac{(\boldsymbol{H}_n\boldsymbol{C}^T)^m}{n^k m!} \prod_{i=1}^{n-\mathrm{m}-1} F_{i+1}(\boldsymbol{C}') \\
&\leq \frac{(\boldsymbol{H}_n\boldsymbol{C}^T)^m (n-m)^k}{n^k m!} \qquad \qquad formula\ (2) \\
&\leq \frac{(\boldsymbol{H}_n\boldsymbol{C}^T)^m}{m!}
\end{aligned}
$$

and

$$
\begin{aligned}
\frac{P\big(O_{\boldsymbol{C}max}(n, m + c_k)\big)}{P\big(O_{\boldsymbol{C}max}(n, m + c_k - 1)\big)} &= \frac{\boldsymbol{H}_n\boldsymbol{C}^T}{mF_{n-m+1}(\boldsymbol{C}')} \\
&\leq \frac{\boldsymbol{H}_n\boldsymbol{C}^T}{m} \qquad \qquad formula\ (3)
\end{aligned}
$$

Since $h_\beta = \binom{k}{\beta}\sum_{j=2}^{n}\frac{1}{(j-1)^\beta}$, we can get $h_0 = n-1$, $h_1 = k\sum_{j=2}^{n}\frac{1}{j-1} \leq k(log(n-1)+\gamma)$ and for $\beta \geq 2$, $h_\beta \leq \binom{k}{\beta}\sum_{j=2}^{n}\frac{1}{(j-1)^2} \leq \binom{k}{\beta}\pi^2/6$. We know $m! > \sqrt{2\pi m}\left(\frac{m}{e}\right)^m$ ($Stirling's\ approximation$).

$$
\frac{P\big(O_{\boldsymbol{C}max}(n, m + c_k)\big)}{P\big(O_{\boldsymbol{C}max}(n, m + c_k - 1)\big)} \leq \frac{c_1 k(log(n-1)+\gamma) + \frac{\pi^2}{6}\sum_{\beta=2}^{k} c_\beta \binom{k}{\beta}}{m}
$$

For $m \geq \left\lceil ekc_1(\log(n-1)+\gamma) + \frac{e\pi^2}{6}\sum_{\beta=2}^{k} c_\beta \binom{k}{\beta}\right\rceil$,

$$
\frac{P\big(O_{\boldsymbol{C}max}(n, m + c_k)\big)}{P\big(O_{\boldsymbol{C}max}(n, m + c_k - 1)\big)} \leq 1/\mathrm{e}
$$

and

$$
\begin{aligned}
P\big(O_{\boldsymbol{C}max}(n, m + c_k)\big) &\leq \frac{(\boldsymbol{H}_n\boldsymbol{C}^T)^m}{m!} \\
&\leq \frac{\left(c_1 k(\log(n-1)+\gamma) + \frac{\pi^2}{6}\sum_{\beta=2}^{k} c_\beta \binom{k}{\beta}\right)^m}{\sqrt{2\pi m}\left(\frac{m}{e}\right)^m} \\
&= \left(\frac{\mathrm{ek}c_1(\log(n-1)+\gamma) + \frac{e\pi^2}{6}\sum_{\beta=2}^{k} c_\beta \binom{k}{\beta}}{m}\right)^m \Big/ \sqrt{2\pi m} \\
&\leq 1/\mathrm{e}
\end{aligned}
$$

For $M = \left\lceil ekc_1(\log(n-1)+\gamma) + \frac{e\pi^2}{6}\sum_{\beta=2}^{k} c_\beta \binom{k}{\beta}\right\rceil + M_1$,

$$
P\big(O_{\boldsymbol{C}max}(n, M + c_k - 1)\big) \leq \exp(-M_1)
$$

So,

$$
\begin{aligned}
P\big(O_{\boldsymbol{C}}(n, m > M)\big) &\le P\big(O_{\boldsymbol{C}max}(n, m > M)\big) \\
&= \sum_{i=M+1}^{n} P\big(O_{\boldsymbol{C}max}(n, i + c_k - 1)\big) \\
&\le \exp(-M_1)
\end{aligned}
$$

For $c_0 = 1$, we can get $c'_0 = 1 - c_0 = 0$ and $P\big(O_{\boldsymbol{C}'}(n, m > M')\big) \le exp(-M_1)$. Hence,

$$
\begin{aligned}
P\big(O_{\boldsymbol{C}}(n, m < n + 1 - M')\big) &= P\big(O_{\boldsymbol{C}}(n, m + c_k - 1 < n - M' + 1 + c_k - 1)\big) \\
&= P\big(O_{\boldsymbol{C}'}(n, m + c'_k - 1 > M' + c'_k - 1)\big) \qquad \textit{Lemma } \mathbf{3.1} \\
&\le \exp(-M_1)
\end{aligned}
$$

**Lemma 5.2.** *Let* $n \ge 2$, $M_1 \in N^+$, $P\big(s_u(n, m)\big) = s_u(n, m)/n!$ *be the probability of* $s_u(n, m)$ *and* $M = \lceil e(log(n-1) + \gamma)\rceil + M_1$, *we get*

$$
P\big(s_u(n, m > M)\big) \le exp(-M_1)
$$

**Proof of Lemma 5.2.** In the proof of Theorem 5.2, we can prove Lemma 5.2 by replacing

$$
\boldsymbol{H}_n \boldsymbol{C}^T \le k c_1(\log(n-1) + \gamma) + \frac{\pi^2}{6} \sum_{\beta=2}^{k} c_\beta \binom{k}{\beta}
$$

with

$$
H(n-1) \le \log(n-1) + \gamma
$$

and $m \ge \left\lceil ekc_1(\log(n-1) + \gamma) + \frac{e\pi^2}{6} \sum_{\beta=2}^{k} c_\beta \binom{k}{\beta} \right\rceil$ with $m \ge \lceil e(\log(n-1) + \gamma)\rceil$.

The following discussion is based on a large $n$. For $c_0 = 0$, the $\boldsymbol{C}$ sequential optimization numbers are nearly concentrated in the head. For $c_0 = 1$, the $\boldsymbol{C}$ sequential optimization numbers are nearly concentrated in the end. For $M_1 \in N^+$, $\exp(-M_1)$ approaches 0 very quickly as $M_1$ increases. Hence, the concentration is noticeable.

Particularly, the probability that the unsigned Stirling numbers of first kind are concentrated in $O(\log n)$ is nearly 100%. This indicates that, in Definition 2.2, for $\boldsymbol{C} = (0,1)$, the probability that the cardinalities of optimization sets are $O(\log n)$ is nearly 100%, which can be useful in the studies using the strategy in Definition 2.1.

**Theorem 5.3.** *For* $k \ge 1$, $n \ge 2$ *and* $0 \le m \le n$, $O_{\boldsymbol{C}}(n, m) \le O_{\boldsymbol{C}'max}(n, n - m)$,

$$
\frac{\sum_{m=0}^{n} O_{\boldsymbol{C}max}(n, m)}{\sum_{m=0}^{n} O_{\boldsymbol{C}}(n, m)} \le exp(\lambda)
$$

*and*

$$
\frac{\sum_{m=0}^{n} O_{\boldsymbol{C}'max}(n, m)}{\sum_{m=0}^{n} O_{\boldsymbol{C}}(n, m)} \le exp(\lambda')
$$

*where* $\lambda = c_0(n-1) + c_1 k(log(n-1) + \gamma) + \frac{\pi^2}{6} \sum_{\beta=2}^{k} c_\beta \binom{k}{\beta}$ *and* $\lambda' = c'_0(n-1) + c'_1 k(log(n-1) + \gamma) + \frac{\pi^2}{6} \sum_{\beta=2}^{k} c'_\beta \binom{k}{\beta}$.

*In particular, for* $s_u(n, m) = O_{(0,1)}(n, m)$,

$$\frac{\sum_{m=0}^{n} s_{umax}(n,m)}{\sum_{m=0}^{n} s_u(n,m)} \le \frac{(n-1)}{n} exp\left(\sum_{j=2}^{n} \frac{1}{j-1} - log(n-1)\right) \le exp(\gamma) \le 1.7811$$

**Proof of Theorem 5.3.**

$$O_{\boldsymbol{C}}(n,m) = O_{\boldsymbol{C}'}(n,n-m) \le O_{\boldsymbol{C}'max}(n,n-m) \qquad Lemma\ 3.1$$

Let $\lambda = c_0(n-1) + c_1 k(\log(n-1)+\gamma) + \frac{\pi^2}{6}\sum_{\beta=2}^{k} c_\beta \binom{k}{\beta}$ and $\lambda' = c_0'(n-1) + c_1' k(\log(n-1) + \gamma) + \frac{\pi^2}{6}\sum_{\beta=2}^{k} c_\beta' \binom{k}{\beta}$. We can get $\boldsymbol{H}_n\boldsymbol{C}^T \le \lambda$ and $\boldsymbol{H}_n\boldsymbol{C'}^T \le \lambda'$. Then,

$$\begin{aligned}
\frac{\sum_{m=0}^{n} O_{\boldsymbol{C}max}(n,m)}{\sum_{m=0}^{n} O_{\boldsymbol{C}}(n,m)} &= \frac{1}{n!^k}\sum_{m=1}^{n}\frac{(n-1)!^k}{(m-1)!}(\boldsymbol{H}_n\boldsymbol{C}^T)^{m-1}\prod_{i=1}^{n-\mathrm{m}} F_{i+1}(\boldsymbol{C}') \\
&\le \frac{1}{n^k}\sum_{m=1}^{n}\frac{\lambda^{i-1}}{(m-1)!}(n-m+1)^k \qquad formula\ (2) \\
&\le \sum_{m=1}^{n}\frac{\lambda^{i-1}}{(m-1)!} \\
&\le \exp(\lambda) \qquad (Taylor\ theorem)
\end{aligned}$$

and

$$\begin{aligned}
\frac{\sum_{m=0}^{n} O_{\boldsymbol{C}'max}(n,m)}{\sum_{m=0}^{n} O_{\boldsymbol{C}}(n,m)} &= \frac{1}{n!^k}\sum_{m=1}^{n}\frac{(n-1)!^k}{(m-1)!}(\boldsymbol{H}_n\boldsymbol{C}'^T)^{m-1}\prod_{i=1}^{n-\mathrm{m}} F_{i+1}(\boldsymbol{C}) \\
&\le \exp(\lambda')
\end{aligned}$$

Particularly, for $s_u(n,m) = O_{(0,1)}(n,m)$,

$$\begin{aligned}
\frac{\sum_{m=0}^{n} s_{umax}(n,m)}{\sum_{m=0}^{n} s_u(n,m)} &= \frac{1}{n!}\sum_{m=1}^{n}\frac{(n-1)!}{(m-1)!}(\boldsymbol{H}_n\boldsymbol{C}^T)^{m-1}\prod_{i=1}^{n-\mathrm{m}} F_{i+1}(\boldsymbol{C}') \\
&\le \frac{1}{n}\sum_{m=1}^{n}\frac{\lambda^{i-1}}{(m-1)!} \\
&\le \exp(\lambda)/n \qquad (Taylor\ theorem) \\
&= \frac{1}{n}\exp\left(\sum_{j=2}^{n}\frac{1}{j-1}\right)
\end{aligned}$$

Let

$$a_n = \frac{1}{n}\exp\left(\sum_{j=2}^{n}\frac{1}{j-1}\right)$$

and

$$f(x) = \frac{x\exp(1/x)}{x+1}$$

where $n \in N^+$ and $x > 0$.

$$\frac{a_{n+1}}{a_n} = \frac{n\exp(1/n)}{n+1}$$

$$f'(x) = -\frac{\exp(1/x)}{x(x+1)^2} < 0$$

$$\lim_{x\to+\infty} f(x) = 1$$

Hence, $f(x) > 1, \frac{a_{n+1}}{a_n} > 1$ and

$$\begin{aligned}\frac{\sum_{m=0}^{n} s_{umax}(n,m)}{\sum_{m=0}^{n} s_u(n,m)} &\le \lim_{n\to+\infty} \frac{1}{n}\exp\left(\sum_{j=2}^{n}\frac{1}{j-1}\right) \\ &= \lim_{n\to+\infty} \frac{n-1}{n}\exp\left(\sum_{j=2}^{n}\frac{1}{j-1} - \log(n-1)\right) \\ &= \exp(\gamma) \\ &\le 1.7811\end{aligned}$$

## 6. Conclusion

In this work, we give the more basic definition of the Stirling numbers of first kind, which yielded the discovery of $\boldsymbol{C}$ sequential optimization numbers. This more basic definition has the potential to discover or prove several more sequences or properties. We deal with multiple inequalities succinctly and roughly, and some of them still have much room for improvement. The unsigned Stirling numbers of first kind, as a special case of $\boldsymbol{C}$ sequential optimization numbers, have been studied in multiple works. We conclude that $\boldsymbol{C}$ sequential optimization numbers have the significance to be studied extensively. For all the $k$+1-tuple vector $\boldsymbol{C}$, there are $2^{k+1}$ kinds of sequences and some of them can have special properties such as $O_{(0,1,1,\dots,1)}(n,m)$ in Lemma 2.1.

## Acknowledgments

The author would like to thank the referees for their careful reading of the manuscript and for their comments and suggestions.